\documentclass[letterpaper, oneside, 12pt]{article}
\usepackage[round]{natbib}

\RequirePackage[OT1]{fontenc}
\RequirePackage{amsthm,amsmath}
\PassOptionsToPackage{hyphens}{url}\usepackage[hidelinks]{hyperref}

\usepackage{mathrsfs}
\usepackage{appendix}
\usepackage{subfigure}
\usepackage{float}
\usepackage{graphicx}
\usepackage{caption}
\bibpunct{(}{)}{;}{a}{}{,}

\usepackage{array,multirow}

\usepackage{amsmath}
\usepackage{amssymb}
\usepackage{graphicx}
\usepackage{amsfonts}
\usepackage{latexsym}
\usepackage{subfigure}
\usepackage{url}
\usepackage{algorithmic}\label{•}
\usepackage{longtable}
\usepackage{multirow}
\usepackage{fancyvrb}
\usepackage{float}
\usepackage{setspace}
\usepackage[hmargin=1.0in,vmargin=1.12in]{geometry}
\usepackage{setspace}
\usepackage{rotating}
\usepackage{enumerate}
\usepackage{comment}
\numberwithin{equation}{section}
\theoremstyle{plain}

\newtheorem*{theorem*}{Theorem}

\newtheorem*{thmdef*}{Theorem}

\newtheorem*{defn*}{Definition}
\newtheorem*{prop*}{Proposition}

\newcommand{\calS}{{\cal S}}

\newcommand{\eps}{\epsilon}
\newcommand{\sigs}{\sigma^2}

\newcommand{\given}{\,|\,}

\newcommand{\sig}{\sigma}

\begin{document}

\title{Bayesian Inference for High Dimensional Changing Linear Regression with Application to Minnesota House Price Index Data}
\author{Abhirup Datta, Hui Zou and Sudipto Banerjee 
}
\date{}
\maketitle

\begin{abstract}
In many applications, the dataset under investigation exhibits heterogeneous regimes that are more appropriately modeled using piece-wise linear models for each of the data segments separated by change-points. Although there have been much work on change point linear regression for the low dimensional case, high-dimensional change point regression is severely underdeveloped.
Motivated by the analysis of Minnesota House Price Index data,  we propose a fully Bayesian framework for fitting changing linear regression models in high-dimensional settings. Using segment-specific shrinkage and diffusion priors, we deliver full posterior inference for the change points and simultaneously obtain posterior probabilities of variable selection in each segment via an efficient Gibbs sampler. Additionally, our method can detect an unknown number of change points and accommodate different variable selection constraints like grouping or partial selection. We substantiate the accuracy of our method using simulation experiments for a wide range of scenarios. We apply our approach for a macro-economic analysis of Minnesota house price index data. The results strongly favor the change point model over a homogeneous (no change point) high-dimensional regression model. 
\end{abstract}


\section{Introduction}\label{Sec: Intro} Modern statistical modeling and inference continue to evolve and be molded by the emergence of complex datasets, where the dimension of each observation in a dataset substantially exceeds the size of the dataset. Largely due to recent advances in technology, such high dimensional datasets are now ubiquitous in fields as diverse as genetics, economics, neuroscience, public health, imaging, and so on. One important objective of high dimensional data analysis is to segregate a small set of regressors, associated with the response of interest, from the large number of redundant ones. Penalized least square approaches like Lasso \citep{lasso}, SCAD \citep{scad}, Elastic Net \citep{enet}, adaptive Lasso\citep{adalasso} etc. are widely employed for high dimensional regression analysis. Bayesian alternatives typically proceed by using hierarchical priors for the regression coefficients aimed at achieving variable selection. Bayesian variable selection methods include stochastic search variable selection \citep{svss}, spike and slab prior \citep{ishwar05}, Bayesian Lasso \citep{blasso}, horseshoe prior \citep{horse}, shrinkage and diffusion prior \citep{naveen14} among others.

Most of the aforementioned approaches assume a single underlying model from which the data is generated. Such homogeneity assumptions may be violated in systems where the variables involved exhibit dynamic behavior and interactions. Common examples include economic time series \citep{gupta97, kezim04, lenard06}, climate change data \citep{climate07}, DNA micro-array data \citep{veera10} and so on. Change point models provide a convenient depiction of such complex relationships by splitting the data based on a threshold variable and using a homogeneous model for each segment. There exists enormous literature on Bayesian methodology addressing various change point problems \citep[see for example][among others]{brad92, barry93, mc93, adams07, turner09}.

Changing linear regression models are a subclass of change point problems, where the linear model relating the response to the predictors varies over different segments of the data. 
Segmentation of the dataset is typically based on unknown change points of a threshold variable like time or age or some other contextual variable observed along with the data. Economic datasets constitute a major domain of application of changing linear models. Many economic time series datasets may be collected over different political and financial regimes, thereby containing several change points with respect to the association with the predictors. In a low dimensional setting, \cite{brad92} used Gibbs' sampling techniques for changing linear models to deliver fully Bayesian inference about the location of the change points and the regression coefficients for each segment. When the set of possible predictors is high dimensional, an additional objective is to identify the (possibly different) sparse supports for each segment. Despite the abundance of Bayesian literature on high dimensional regression and on change point models, it appears there is no extant Bayesian work on high dimensional changing linear regression.

This manuscript intends to bridge this gap by proposing a hierarchical methodology for high dimensional changing linear models. We embed Bayesian variable selection techniques in a change point setup to  simultaneously detect the number and location of the change points as well as to identify the true sparse support for each of the linear models.  We use teh newly proposed shrinkage and diffusion priors \citep{naveen14} for the regression coefficients in each segment to perform variable selection. We provide an efficient Gibbs' sampler that delivers full posterior inference on the change points, posterior selection probabilities for each variable for all segments and posterior predictive distributions for the response. Our fully Bayesian approach is flexible to the choice of variable selection priors and offers the scope for several structural modifications tailored to specific data applications. For example, constraints like grouping the selection of a variable across all the segments can be easily achieved using group selection priors. Other constraints like partial selection within or between the segments can also be accommodated in our setup. Numerical studies reveal that for a wide range of scenarios, our proposed methodology can accurately detect the change points and select the correct set of predictors. We demonstrate the applicability of our method for a macro-economic analysis of Minnesota house price index data. The results strongly favor our change point model over a homogeneous high dimensional regression model.

Classical penalized least square approaches mentioned earlier can also be used in a change point setup. By treating the unknown change points as additional tuning parameters, one can split the data using fixed values of these change points and use some penalized loss function to achieve variable selection for each segment. For example, \cite{cplasso} uses Lasso penalty to estimate the coefficients for each segment. Subsequent application of cross validation or model selection techniques will yield the optimal change points from a grid of possible values. However, our fully Bayesian approach has several advantages over this. Firstly, the grid search approach is computationally highly inefficient especially for more than one change points. On the other hand, a prior specification for the change points in our Bayesian model enables standard MCMC techniques to efficiently generate posterior samples. Moreover, in many real applications, change in association between variables can occur over a range of the threshold variable. Point estimates of change points obtained from classical approaches fail to accurately depict such scenarios. Bayesian credible intervals obtained from the posterior distributions provide a much more realistic quantification of the uncertainty associated with the location of the change points, which is very difficulty to accomplish if one uses the grid search approach.

The rest of the manuscript is organized as follows. In Section \ref{sec:method} we present our methodology in details including extensions to unknown number of change points and alternate prior choices. Results from several simulated numerical studies are provided in Section \ref{sec:num}. In Section \ref{sec:minn} we present the details of a house price index data analysis using our change point methodology. We conclude this paper with a brief review and pointers to future research.

\section{Method}\label{sec:method}
We consider a traditional high dimensional setup with the $n\times 1$ response vector $y=(y_1,y_2,\ldots,y_n)'$ and corresponding $n\times p$ covariate matrix $X=(x_1,x_2,\ldots,x_n)'$ where $p$ can be larger than $n$. We further assume that for every observation $y_i$ we observe another quantitative variable $t_i$ such that the association between $y_i$ and $x_i$ depends on the values of $t_i$. In a linear regression setup, this dynamic relationship between the response $y_i$ and the corresponding $p\times 1$ vector of covariates $x_i$ can be expressed as $E(y_i \given x_i, t_i) =x_i'\beta_k$ for all $i$ such that $\tau_{k-1} < t_i < \tau_k$ where $\tau_0 < \tau_1 < \ldots < \tau_K < \tau_{K+1}=n$. The change-points $\tau_1, \tau_2, \ldots, \tau_K$ are typically unknown while the number of change-points $K$ may or may not be known depending on the application.

As the number of regressors ($p$) is large, our goal is to select the relevant variables for this regression. However, for this changing linear regression, the set of relevant regressors may depend on the value of the threshold variable $t$ and variable selection procedures applied disregarding the dependence on $t$ can lead to erroneous variable selection. Let $S_k$ denotes the support of $\beta_k$ where $s_k=|S_k|$ is typically much less than $p$. Our goal is to simultaneously detect the change-points $\tau_k$ and estimate $S_k$ for all $k=1,2,\ldots,K$. We initially assume only one change-point $\tau$ i.e. $K=1$. Extensions to more than one (and possibly unknown number of) change points are discussed later in Section \ref{sec:mult}.

\subsection{One Change Point Model}\label{sec:cp1} We assume a changing linear regression model
\begin{align}\label{eq:model}
y_i= \left \{ \begin{array}{c}
 x_i'\beta_1 + \eps_i \mbox{ if } t_i \leq \tau \\
 x_i'\beta_2 + \eps_i \mbox{ if } t_i > \tau
\end{array} \right.
\end{align}
where $\beta_1$, $\beta_2$ are both sparse $p\times 1$ vectors such that  $\beta_1 \neq \beta_2$ and $\eps_i \sim N(0,\sigs)$  denotes the independent and identically distributed noise. In order to accomplish variable selection both before and after the change point, we use shrinking and diffusion (BASAD) priors proposed in \cite{naveen14} for  $\beta_1$ and $\beta_2$. To be specific, we assume $\beta_k \given Z_k, \sigs \sim N(0,\sigs diag(\gamma_{1k} Z_k  + \gamma_{0k}(1-Z_k) ))$ for $k=1,2$ where $Z_k=(Z_{k1},Z_{k2},\ldots,Z_{kp})'$  is a $p\times 1 $ vector of zeros and ones. The hyper-parameters $\gamma_{0k}$ and $\gamma_{1k}$ are scalars chosen to be very small and very large respectively. Hence, $\beta_{kj}$---the $j^{th}$ component of $\beta_k$---is assigned a shrinking prior if $Z_{kj}$ equals $0$ and a diffusion (flat) prior if $Z_{kj} =1$. $Z_{kj}$'s are assumed to be apriori independent each following $Bernoulli(q_k)$. Hence $q_k$ controls the prior model size for the $k^{th}$ segment. The choices for the hyper-parameters $\gamma_{0k}$, $\gamma_{1k}$ and $q_k$ are discussed in Section~\ref{sec:num}. We assume a uniform prior for the change-point $\tau$ and a conjugate Inverse Gamma prior for the noise variance $\sigs$. The full Bayesian model can now be written as:
\begin{align}\label{eq:like}
\prod_{i: t_i \leq \tau} & N( y_i  \given x_i'\beta_1, \sigs) \prod_{i: t_i > \tau} N(y_i \given x_i'\beta_2, \sigs) \times Unif(\tau \given a_\tau, b_\tau) \times IG(\sigs \given a_\sigma, b_\sigma) \times \nonumber \\
& \prod_{k=1}^2 \left( N(\beta_k \given 0,\sigs diag(\gamma_{1k} Z_k  + \gamma_{0k}(1-Z_k) )) \times \prod_{j=1}^p Bernoulli(Z_{kj} \given q_k) \right)
\end{align}

We use Gibbs' sampler to obtain posterior samples of $\tau$, $\beta_k$ and $Z_k$ for $k=1,2$. Let $\tau \given \cdot$ denote the full-conditional distribution of $\tau$ in the Gibbs' sampler. We use similar notation to denote the other full conditionals. Let $U_1=\{i \given t_i \leq \tau \}$ and $U_2=\{ i \given t_i > \tau\}$. For $k=1,2$, let $Y_k$ and $X_k$ denote the response vector and covariate matrix obtained by stacking up the observations $U_k$. From the full likelihood in (\ref{eq:like}), we have
\begin{align*}
\beta_k \given \cdot \sim  N(V_k X_k' Y_k , \;  \sigs V_k) & \mbox{ where } V_k = (X_k'X_k+diag(\gamma_{1k} Z_k  + \gamma_{0k}(1-Z_k) )^{-1} )^{-1} \\
\sigs  \given \cdot \sim  IG&(a_\sigma+n/2, b_\sig + +\frac 12 \sum_{k=1}^2 ||Y_k-X_k\beta_k||^2) \\
p(\tau \given \cdot) \propto \prod_{k=1}^2 & \prod_{i \in U_k} N(y_i  \given x_i'\beta_k, \sigs) \times Unif(\tau \given a_\sigma, b_\sigma)\\
Z_{kj} \given \cdot \sim Bernoulli & \left(\frac {q_k \; \phi(\beta_{kj}/\sqrt{\sigs\gamma_{1k}})}{q_k\;\phi(\beta_{kj}/\sqrt{\sigs\gamma_{1k}})+(1-q_k)\;\phi(\beta_{kj}/\sqrt{\sigs\gamma_{0k}})} \right)
\end{align*}
where $\phi(x)$ denotes the density of standard normal distribution. We observe that the full conditionals of $\beta_k$, $Z_{kj}$ and $\sigs$ follow conjugate distributions and are easily updated via the Gibbs' sampler. Only $p(\tau \given \cdot)$ does not correspond to any standard likelihood and we use a Random Walk Metropolis-Hastings step within the Gibbs' sampler to update $\tau$.

\subsection{Multiple change points}\label{sec:mult} So far we have limited our discussion to the presence of only one change point. However, our method can be easily extended to multiple change points. If we have $K$ change points $\tau_1 < \ldots < \tau_K $, the joint likelihood in (\ref{eq:like}) can be generalized to
\begin{align}
\prod_{k=1}^K \left( \prod_{i:\tau_{k-1} < t_i \leq \tau_k}  \right. N(y_i \given x_i'\beta_k, \sigs) \times N(\beta_k \given 0,\sigs diag(\gamma_{1k} Z_k &+ \gamma_{0k}(1-Z_k) )) \times  \nonumber  \\
 \left.  \prod_{j=1}^p Bernoulli(Z_{kj} \given q_k) \right) \times p(\tau_1,\tau_2,\ldots,\tau_K) \times & IG(\sigs \given a_\sigma, b_\sigma)
\end{align}
To ensure identifiability of the change points, the prior $p(\tau_1,\tau_2,\ldots,\tau_K)$ should be supported on $\tau_1 < \tau_2 < \ldots < \tau_K$. To accomplish this we choose $p(\tau_1,\tau_2,\ldots,\tau_K)$ as the density of ordered statistics of a sample of size $K$ from $Unif(a_\tau,b_\tau)$. The Gibbs' sampler remains essentially same as in Section \ref{sec:method} with the Metropolis random walk step now being used to update the entire change point vector $(\tau_1,\tau_2,\ldots,\tau_K)'$.

\subsection{Determining the number of change points}

However, often in applications, the number of change points is unknown. In our fully Bayesian approach this can potentially be handled by adding a prior for the number of change points ($K$). Introducing this additional level of hierarchy comes with the caveat that different values of $K$ yields parameter sub-spaces of different sizes and interpretations. To elucidate, a one change point model splits the data into two segments, with separate coefficient vectors $\beta_1$ and $\beta_2$, creating a parameter space of dimension $2p$ whereas a no change point model has a single $\beta$ of dimension $p$ with a possible interpretation that it is some average of $\beta_1$ and $\beta_2$ over the two segments. Therefore, a Markov Chain Monte Carlo sampling for $K$ will involve jumping within and between different sub-spaces.

\cite{green95} proposed the extremely general and powerful reversible jump MCMC (RJMCMC) sampler for sampling across multiple parameter spaces of variable dimensions. We can seamlessly adopt an RJMCMC joint sampler to obtain the posterior distribution for the number of change points. When naively implemented, RJMCMC experiences poor acceptance rates for transitions to parameter sub-spaces with different dimensionality. This leads to widely documented convergence issues \citep{green09,sisson11}. The problem will be exacerbated in our setup due to the high dimensionality of the parameter spaces.

Several improvements and alternatives to RJMCMC have been proposed over the years including efficient proposal strategies to effectuate frequent cross-dimensional jumps \citep{rich97,brooks03,ehlers08,farr15}, product space search \citep{brad95,della02} and parallel tempering \citep{litt09}.
All these approaches can be adapted in our setup to determine the number of change points. However, many of these approaches are accompanied by their own computational burden such as running several chains or apriori obtaining posterior distributions for each individual model before running the joint sampler. We concur with \cite{brad01} and \cite{green12} that it is often expedient to use simpler model selection approaches based on individual models. Hence, popular Bayesian model comparison metrics like DIC \citep{spieg02} and l-measure \citep{gelf98} remains relevant to select the number of change points in our case. For example, if $\theta$ is the complete set of parameters associated with the model, for each $K$ we can compute the DIC score
\begin{equation}
\mbox{DIC} = 2E\left(D ( y \given \theta) \given y\right) - D\left(y \given E(\theta \given y)\right) = E\left(D ( y \given \theta) \given y\right)+p_D
\end{equation}
where $D(y \given \theta)$ is the deviance function 
and $p_D=E\left(D ( y \given \theta) \given y\right) - D\left(y \given E(\theta \given y)\right)$ is interpreted as effective sample size. Hence, DIC penalizes more complex models and is particularly suitable for our change point context where higher number of change points will lead to overfitting.
Parallel computing can be utilized to simultaneously run the MCMC sampler for different values of $K$ and then the optimal $K$ can be selected as the one yielding lowest DIC score.

All the methods for selecting the number of change points discussed here can be used in conjunction with our approach. It is prudent to predicate the choice on the nature of the application at hand and the computational resources available.



\subsection{Alternate prior choices}\label{sec:prior} 
We observe from Equation \ref{eq:like} that conditional on the value of the change point $\tau$, the joint likelihood can be decomposed into individual likelihoods for the regression before and after the change point along with the corresponding priors for the regression coefficients. This allows for a lot of flexibility in the choice of  priors for the regression coefficients.

We have focused on the BASAD prior. One can also use other priors to achieve variable selection. For example, using Laplace (double exponential) priors for the $\beta_k$'s will yield a Bayesian Lasso \citep{blasso} with change point detection. To facilitate the discussion, consider the one change point model. By using the Laplace prior, the full hierarchical specification for the coefficient vectors $\beta_k$ for $k=1,2$ can be specified as:
\begin{align}\label{eq:blasso}
\beta_k \given \sigs, \eta_k & \overset{ind}{\sim} N(0,\sigs \;  diag(\eta_k)) \mbox{ where } \eta_k=(\eta_{k1},\eta_{k2},\ldots,\eta_{kp})' \nonumber \\
\eta_{kj} \given  \lambda_k & \overset{ind}{\sim} Exp(\lambda_k^2/2) \mbox{ and }
\lambda_k^2 \sim Gamma(r_k,s_k)
\end{align}
The prior specification for $\sigs$ and $\tau$ can be kept same as in (\ref{eq:like}). The Gibbs' sampler for the Bayesian Lasso provided in \cite{blasso} can now be used to sample from the following full conditionals:
\begin{align*}
\beta_k \given \cdot \sim N(V_kX_k'y_k,\sigs V_k) &\mbox{ where } V_k=(X_k'X_k+diag(\eta_k)^{-1})^{-1} \\
\sigs \given \cdot \sim IG(a_\sigma+n/2, & \;b_\sigma+\frac 12 \sum_{k=1}^2 ||Y_k-X_k\beta_k||^2)\\
1/\eta_{kj} \given \cdot \sim \mbox{ Inv-} & \mbox{Gauss  } \left(\sqrt \frac{\lambda_k^2\sigs}{\beta_{kj}^2}, \lambda_k^2 \right) \\ \lambda_k^2 \given \cdot \sim Gamma&(r_k+p/2,s_k + \frac 12 ||\eta_k||_2^2) \\
p(\tau \given \cdot) \propto \prod_{k=1}^2 \prod_{i \in U_k} &N(y_i  \given x_i'\beta_k, \sigs) \times Unif(\tau \given a_\sigma, b_\sigma)
\end{align*}

Additional information regarding grouping or structuring of the variables are often available in the context of variable selection. In presence of a change point, additional constraints can specify grouped selection both within and/or between the $\beta_k$'s. For example, in a single change point setup, it may be plausible that the set of relevant variables remain unchanged before and after  the change point, with change occurring only with respect to the strength of association between $y_i$ and $x_i$. Such additional structural constraints both within and across $\beta_k$'s can easily be accommodated in our setup via a suitable choice of prior. To elucidate, we can rewrite (\ref{eq:model}) as $y_i=z_i(\tau)'\zeta + \eps_i$  where $z_i=(I(t_i \leq \tau)  x_i', I(t_i > \tau) x_i')' $ and $\zeta=(\beta_1', \beta_2')'$. To incorporate the constraint that $\beta_1$ ans $\beta_2$ share the same support, one can use a Bayesian group lasso \citep{bgrplasso} with M-Laplace priors on the groups $\zeta_j=(\beta_{1j},\beta_{2j})'$ for $j=1,2,\ldots,p$. The M-Laplace prior \[ p(\zeta_j \given \sigs, \lambda^2) \propto \frac {2\lambda^2}{\sigs} \exp(-\sqrt \frac {2\lambda^2}{\sigs} ||\zeta_j||_2)\] has a convenient two-step hierarchical specification:
\begin{equation}\label{eq:grplasso}
\zeta_j \given \eta_j \overset{ind}{\sim} N(0, \sigs \eta_j I) ;\; \eta_j \given \lambda^2 \overset{ind}{\sim} Gamma(3/2, \lambda^2) ; \; \lambda^2 \sim Gamma(r,s)
\end{equation}
The full conditional distributions of the parameters provided in \cite{bgrplasso} can now be used to implement the Gibbs' sampler with the additional Metropolis Random walk step for updating the change point $\tau$. Any other information like hierarchical selection or anti-hierarchical selection both within and between the $\beta_k$'s can also be accommodated via suitable priors.

Often, in real data applications, prior knowledge dictates the inclusion of certain variables in the model and variable selection is sought only for the remaining variables. 
Such constraints can be easily achieved in our setup by using standard Gaussian prior for that specified subset and BASAD prior for the remaining variables.

\subsection{Variable selection after MCMC}
When there are finite many candidate models, Bayesian model selection typically proceeds by selecting the candidate model with the highest posterior probability.
However, in our setup the regression coefficients are continuous.
For variable selection, we use the median probability model \citep{medianprob} which is computationally easy and is optimal in terms of prediction. To be specific,
$\beta_{kj}$ is included in the model if the posterior probability of $Z_{kj}=1$ is greater than $0.5$.

\section{Numerical Studies}\label{sec:num} We conducted numerical experiments to illustrate the performance of our method both for single and multiple change points. For all the simulation studies we used  $100,000$ MCMC iterations. Multiple chains were run with different choices of initial values and convergence was typically achieved within the first $20,000$ iterations. Nevertheless, we discarded the first $50,000$ as burn-in and used the subsequent $50,000$ samples for inference.

\subsection{One change point}\label{sec:uni} We assume $t_i=i$ and generate data from the model $y_i=N(x_i'\beta_1,\sigs)$ for $ i \leq \tau$ and $y_i=N(x_i'\beta_2,\sigs)$ for $i > \tau$ where $\beta_1=(3,1.5,0,0,2,0,\ldots,0)$ and $\beta_2=-\beta_1$. The rows of $X$ were independent and identically distributed normal random  variables with zero mean and covariance $\Sigma_X$. Two structures were used for $\Sigma_X$ --- auto-regressive (AR) with $\Sigma_{X,ij}=0.5^{|i-j|}$ and compound symmetry (CS) with $\Sigma_{X,ij}=0.5+0.5 I(i=j)$. The noise variance $\sigs$ was fixed at $1$ and the sample size was chosen to be $200$. Two different model sizes --- $p=250$ and $p=500$ were used. The change point $\tau$ was chosen to vary between $50.5$, $100.5$ and $150.5$. Since, the sample size is $200$, these three choices of $\tau$ respectively correspond to changes in the regression model at the initial , middle or later portion of the data. We used three different prior choices for the coefficients --- the BASAD prior, the Bayesian Lasso prior \citep{blasso} and the Bayesian Group Lasso prior \citep{bgrplasso}. The last choice was used to investigate any possible benefits of using a grouped variable selection as it is known that $\beta_1$ and $\beta_2$ has same support. The range of the uniform prior for $\tau$ was chosen to be $(20,180)$ and a normal proposal density with tuning variance of $0.1$ was used for the Metropolis update of $\tau$. The prior for $\sigs$ was chosen to be $IG(2,1)$. The hyper-parameters $\gamma_{0k}$, $\gamma_{1k}$ and $q_k$ were chosen as follows. Let $\tau_0$ denote the initial estimate for $\tau$. Then $n_1=[\tau_0]$ and $n_2=n-n_1$ denotes the initial sample sizes for the two segments. We used
\[ \gamma_{0k}=\frac{\hat \sigma_k^2} {10n_k} \qquad , \qquad \gamma_{1k}=\hat \sigma^2_k\; \max \left(\frac {p^{2.1}} {100n_k},\log n_k \right) \]
where $\hat \sigma^2_k$ was the sample variance of $Y_k$ for $k=1,2$. The hyper-parameters $q_k$ were chosen such that the prior model sizes $\sum_{j=1}^p Z_{kj}$ were greater than $\min(p-1, \max(10, \log n_k))$ with probability  $0.1$. These choices of $\gamma_{0k}$, $\gamma_{1k}$ and $q_k$ were adapted from \cite{naveen14}.

The posterior median estimates of the change points for all the scenarios are provided in Tables \ref{tab:tau250} ($p=250$) and \ref{tab:tau500} ($p=500$).
\begin{table}[h!]
\centering
\caption{Single change point model with $n=200$ and $p=250$: Posterior median estimates (and $95\%$ confidence intervals) of $\tau$ using BASAD, Bayesian Lasso (BL) and Bayesian Group Lasso (BGL) priors.}\label{tab:tau250}
\begin{tabular}{c|c|c|c|c} \hline
& True $\tau$ & BASAD & BL & BGL \\ \hline
 & 50.5 & 50.5 (50.0, 51.0) & 50.4 (48.9, 51.8) & 50.5 (49.5, 51.5) \\
AR & 100.5 & 100.5 (100.0, 101.0) & 100.5 (100.0, 101.0) & 100.5 (100.0, 101.0) \\
 & 150.5 & 150.5 (150.0, 151.0) & 150.5 (150.0, 151.0) & 150.5 (150.0, 151.0) \\ \hline
 & 50.5 & 50.1 (49.1, 51.9) & 50.0 (49.1, 51.8) & 50.1 (49.1, 51.9)\\
CS & 100.5 & 100.5 (100.0, 101.0) & 100.5 (100.0, 101.0) & 100.5 (100.0, 101.0)\\
 & 150.5 & 150.5 (150.0, 151.0) & 150.2 (148.6, 151.0) & 150.4 (149.1, 151.0)\\ \hline
\end{tabular}
\end{table}
\begin{table}[h]
\centering
\caption{Single change point model with $n=200$ and $p=500$: Posterior median estimates (and $95\%$ confidence intervals) of $\tau$ using BASAD, Bayesian Lasso (BL) and Bayesian Group Lasso (BGL) priors.}\label{tab:tau500}
\begin{tabular}{c|c|c|c|c} \hline
& True $\tau$ & BASAD & BL & BGL \\ \hline
 & 50.5 & 50.5 (50.0, 51.0) & 49.7 (47.6, 53.7) & 50.4 (48.2, 53.1)  \\
AR & 100.5 & 101.3 (99.3, 102.0) & 101.3 (99.1, 103.9) & 101.0 (99.1, 103.3) \\
 & 150.5 & 150.5 (150.0, 151.0) & 151.3 (148.7, 152.9) & 150.7 (150.0, 152.8)  \\ \hline
 & 50.5 & 50.3 (49.1, 51.0) & 50.1 (48.5, 51.8) & 50.2 (49.1, 51.5)  \\
CS & 100.5 & 100.3 (99.1, 101.0) & 100.0 (99.1, 100.9) & 99.9 (99.1, 100.9)  \\
 & 150.5 & 150.5 (150.0, 151.0) & 150.5 (150.0, 151.0) & 150.5 (150.0, 151.0) \\ \hline
 \end{tabular}
\end{table}
We observe that all the 3 models estimate the change point with high accuracy.

We then turn our attention to variable selection. Let $C_k$ and $IC_k$ denote the number of true and false regressors respectively selected for the $k^{th}$ segment of the data for $k=1,2$.  As discussed earlier, we used a cut-off of $0.5$ for the posterior probability of the binary $Z_{kj}$'s in the BASAD model to select the variables. The Bayesian Lasso and Group Lasso are devoid of such binary selection parameters and variable selection was based on the posterior confidence intervals i.e $\beta_{kj}$ was not selected if its posterior confidence interval covered zero. Table \ref{tab:cic250} provides the $C_k$ and $IC_k$ numbers for each method for $p=250$.
\begin{table}[h!]
\centering
\caption{Single change point model with $n=200$ and $p=250$: Number of correct and incorrect predictors selected by BASAD, Bayesian Lasso (BL) and Bayesian Group Lasso (BGL). Cases where any method missed at least one true regressor are highlighted using *.}\label{tab:cic250}
\begin{tabular}{c|c|c|ccc|ccc} \hline
&&& \multicolumn{3}{c|}{AR} & \multicolumn{3}{c}{CS}\\ \hline
 &  & True & BASAD & BL & BGL & BASAD & BL & BGL\\ \hline
\multirow{4}{*}{$\tau=50.5$} & $C_1$ & 3 & 3 & 3 & 3  & 3 & 2* & 3\\
 & $C_2$ & 3 & 3 & 3 & 3  & 3 & 3 & 3\\
 & $IC_1$ & 0 & 0 & 0 & 0  & 0 & 0 & 0\\
 & $IC_2$ & 0 & 0 & 2 & 2   & 0 & 2 & 1\\ \hline
\multirow{4}{*}{$\tau=100.5$} & $C_1$ & 3 & 3 & 3 & 3  & 3 & 3 & 3\\
  & $C_2$& 3 & 3 & 3 & 3  & 3 & 3 & 3\\
 & $IC_1$ & 0 & 0 & 0 & 0  & 0 & 0 & 0\\
 & $IC_2$ & 0 & 0 & 0 & 0 & 0 & 0 & 0\\ \hline
\multirow{4}{*}{$\tau=150.5$} & $C_1$ & 3 & 3 & 3 & 3 & 3 & 3 & 3\\
 & $C_2$ & 3 & 3 & 2* & 3 &  3 & 2* & 3\\
 & $IC_1$ & 0 & 0 & 0 & 0  & 0 & 2 & 1\\
 & $IC_2$ & 0 & 0 & 0 & 0 & 0 & 0 & 0\\ \hline
\end{tabular}
\end{table}
 We observe that the BASAD prior achieves perfect variable selection for all the scenarios while the Bayesian LASSO often misses out on a true variable and includes some incorrect predictors. The Bayesian Group LASSO always selects the true set of regressors but often includes one false regressor. For $p=500$, (Table \ref{tab:cic500}), the BASAD once again selects the exact set of predictors.
\begin{table}[!h]
\centering
\caption{Single change point model with $n=200$ and $p=500$: Number of correct and incorrect predictors selected by BASAD, Bayesian Lasso (BL) and Bayesian Group Lasso (BGL). Cases where any method missed at least one true regressor are highlighted using *.}\label{tab:cic500}
\begin{tabular}{c|c|c|ccc|ccc} \hline
&&& \multicolumn{3}{c|}{AR} & \multicolumn{3}{c}{CS}\\ \hline
 &  & True & BASAD & BL & BGL & BASAD & BL & BGL\\ \hline
\multirow{4}{*}{$\tau=50.5$} & $C_1$ & 3 & 3 & 1* & 3 & 3 & 1* & 2* \\
 & $C_2$ & 3 & 3 & 3 & 3 & 3 & 3 & 3 \\
 & $IC_1$ & 0 & 0 & 0 & 0 & 0 & 0 & 0 \\
 & $IC_2$ & 0 & 0 & 0 & 0 & 0 & 0 & 0 \\ \hline
\multirow{4}{*}{$\tau=100.5$} & $C_1$ & 3 & 3 & 3 & 3 & 3 & 3 & 3 \\
 & $C_2$ & 3 & 3 & 3 & 3 & 3 & 3 & 3 \\
 & $IC_1$ & 0 & 0 & 0 & 0 & 0 & 0 & 0 \\
 & $IC_2$ & 0 & 0 & 0 & 0 & 0 & 0 & 0 \\ \hline
\multirow{4}{*}{$\tau=150.5$} & $C_1$ & 3 & 3 & 3 & 3 & 3 & 3 & 3 \\
 & $C_2$ & 3 & 3 & 1* & 3 & 3 & 0* & 3 \\
 & $IC_1$ & 0 & 0 & 0 & 0 & 0 & 0 & 0 \\
 & $IC_2$ & 0 & 0 & 0 & 0 & 0 & 0 & 0 \\ \hline
\end{tabular}
\end{table}
 However, the performance of the Bayesian Group LASSO and especially the Bayesian LASSO worsens with the latter often being able to select only one correct predictor.

\begin{figure}[!h]
\centering
\subfigure[$MSE_1$: AR]{\includegraphics[scale=0.2]{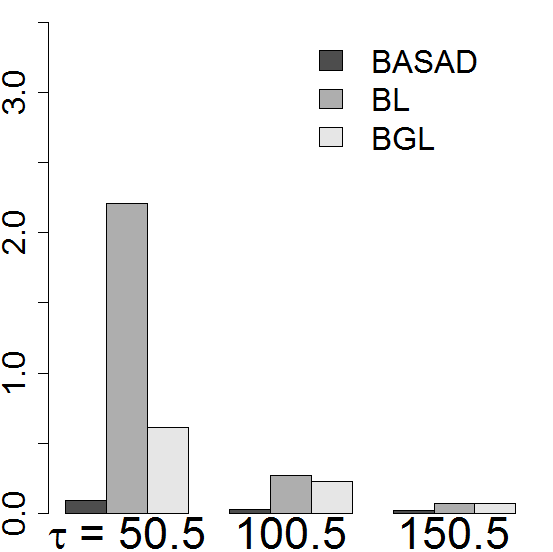}}
\subfigure[$MSE_2$: AR]{\includegraphics[scale=0.2]{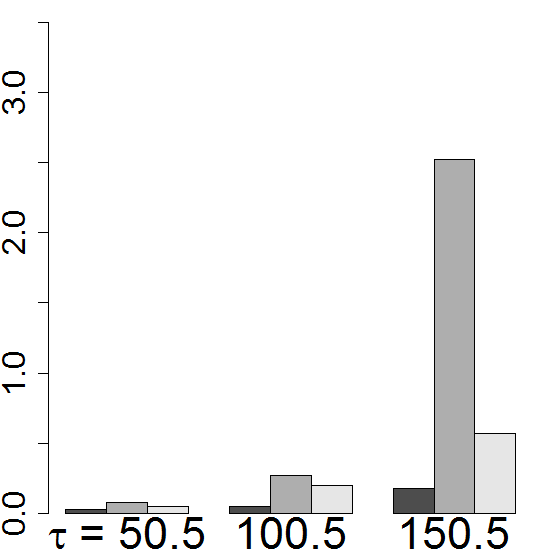}}
\subfigure[$MSE_1$: CS]{\includegraphics[scale=0.2]{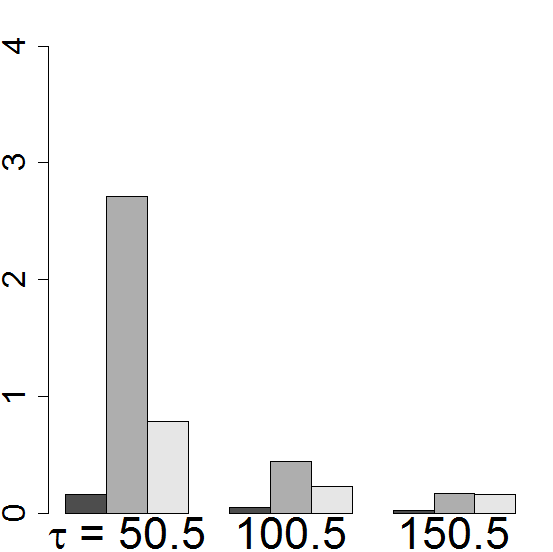}}
\subfigure[$MSE_2$: CS]{\includegraphics[scale=0.2]{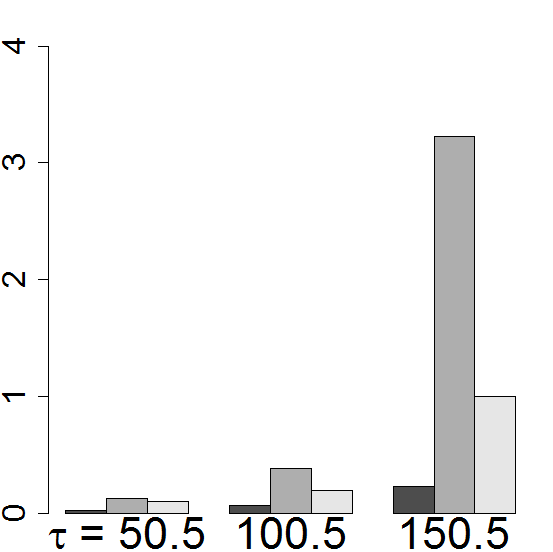}}
\caption{Single change point model with $p=250$: Posterior median estimates of the truncated Mean Squared Error ($MSE_k$) for $\beta_1$ and $\beta_2$ using BASAD, Lasso and Group Lasso (GL) priors}\label{fig:mse1cp250}
\end{figure}

\begin{figure}[!h]
\centering
\subfigure[$MSE_1$: AR]{\includegraphics[scale=0.2]{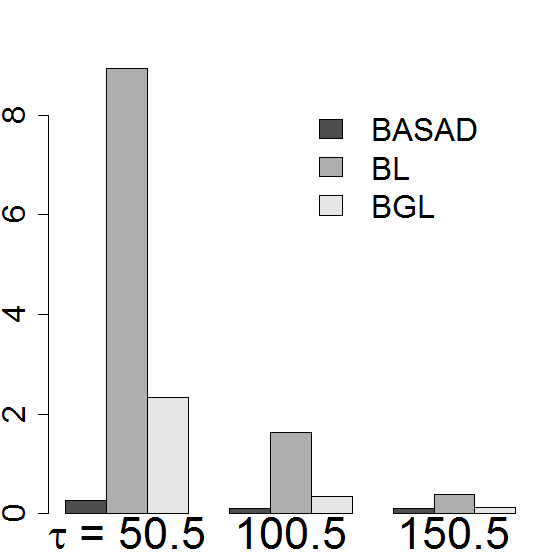}}
\subfigure[$MSE_2$: AR]{\includegraphics[scale=0.2]{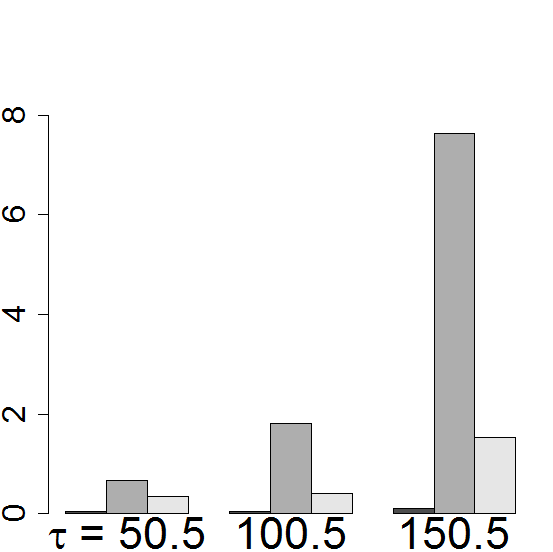}}
\subfigure[$MSE_1$: CS]{\includegraphics[scale=0.2]{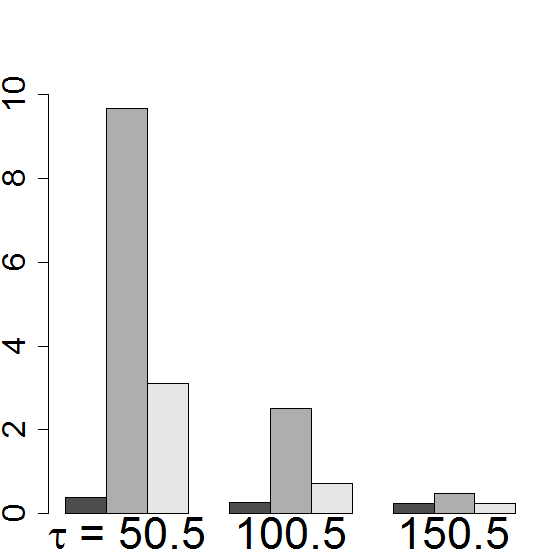}}
\subfigure[$MSE_2$: CS]{\includegraphics[scale=0.2]{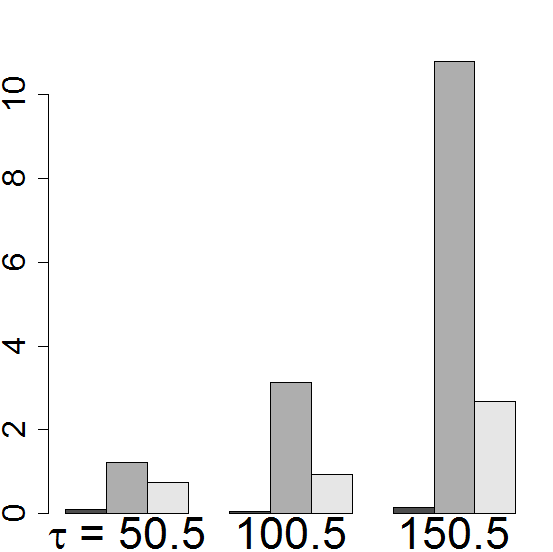}}
\caption{Single change point model with $p=500$ : Posterior median estimates of the truncated Mean Squared Error ($MSE_k$) for $\beta_1$ and $\beta_2$ using BASAD, Lasso and Group Lasso (GL) priors}\label{fig:mse1cp500}
\end{figure}

In additional to the variable selection metrics, we also assess the three methods based on the coefficient estimates for the true predictors using the Mean Squared Error truncated on the true support i.e.
\[ MSE_k = || \beta_k[\calS_k] - \hat\beta_k[\calS_k] ||_2^2 \mbox{ for } k=1,2 \]
 where $\calS_k$ denotes the true support of $\beta_k$ and $\hat \beta_k$ denote its posterior estimate. Figures \ref{fig:mse1cp250} and \ref{fig:mse1cp500} plots the $MSE_k$ numbers for $p=250$ and $p=500$ respectively. $\beta_1[ \calS _1 ]$
We observe that the BASAD stands out with uniformly lowest MSE numbers across all scenarios. It is important to note that, $MSE_1$ tends to be higher when $\tau=50.5$ while $MSE_2$ is higher when $\tau=150.5$. This behavior is expected as for $\tau=50.5$, sample size for estimating $\beta_1$ is effectively $50$ while that for $\beta_2$ is $150$. The Bayesian Lasso seems to be worst impacted by this effective sample size. It performs worse for the Compound Symmetry covariance structure and for higher model size ($p=500$). The Bayesian Group Lasso, enjoying the additional knowledge of the structural constraint, conceivably performs better than the Bayesian Lasso. However, the BASAD seems to be least impacted by effective sample size, model size or covariance structure producing accurate variable selection, change point detection and estimation across all scenarios.

\subsection{Two change points}\label{sec:two} We demonstrate the applicability of our method to multiple change points using a two change point setup. The three coefficient vectors are given by
\begin{align*}
\beta_1=(3,0,0, \ldots ,0)', \; \beta_2=(3,1.5,0,0,\ldots,0)',
\beta_3=(3,1.5,0,0,2,0,0,\ldots,0)'
\end{align*}
Three pairs of values for the change points $(\tau_1,\tau_2)$ are selected --- $(50.5,100.5)$, $(50.5,150.5)$ and $(100.5,150.5)$.
Other specifications including sample size, model size and covariance of the predictors are kept unchanged from Section \ref{sec:uni}. We observed in Figures \ref{fig:mse1cp250} and \ref{fig:mse1cp500}, that the Bayesian Lasso becomes erratic for small effective sample sizes. For 2 change points, the effective sample size is further lowered and the Bayesian Lasso faced convergence issues in the MCMC sampler. The Bayesian Group Lasso also cannot be used here as the coefficient vectors for different segments do not share a common support. Hence, we only present the results for the BASAD prior.

\begin{figure}
\centering
\subfigure[$\tau_1=50.5$, $\tau_2=100.5$, $p=250$]{\includegraphics[scale=0.3]{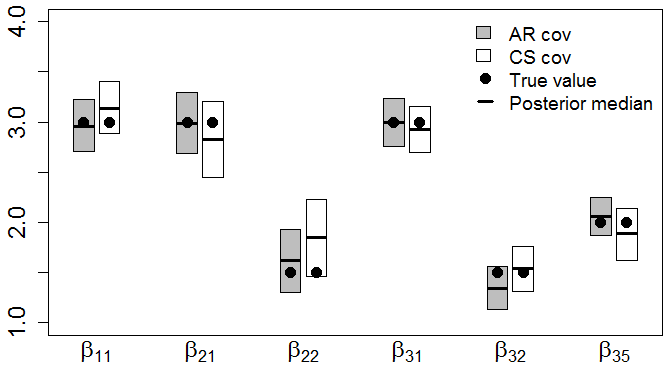}}
\subfigure[$\tau_1=50.5$, $\tau_2=100.5$, $p=500$]{\includegraphics[scale=0.3]{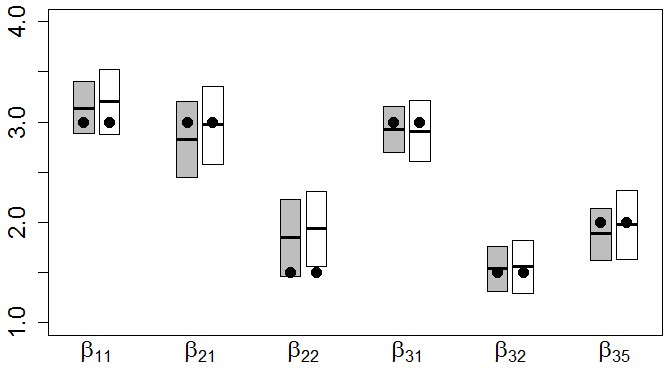}}\\
\subfigure[$\tau_1=50.5$, $\tau_2=150.5$, $p=250$]{\includegraphics[scale=0.3]{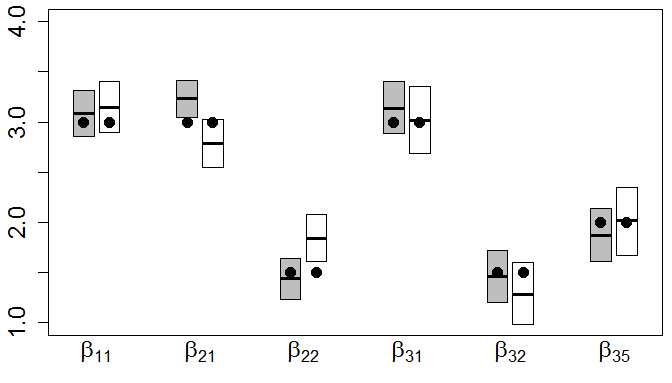}}
\subfigure[$\tau_1=50.5$, $\tau_2=150.5$, $p=500$]{\includegraphics[scale=0.3]{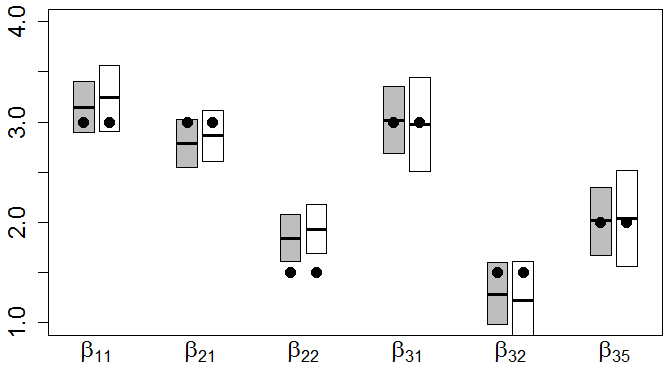}}\\
\subfigure[$\tau_1=100.5$, $\tau_2=150.5$, $p=250$]{\includegraphics[scale=0.3]{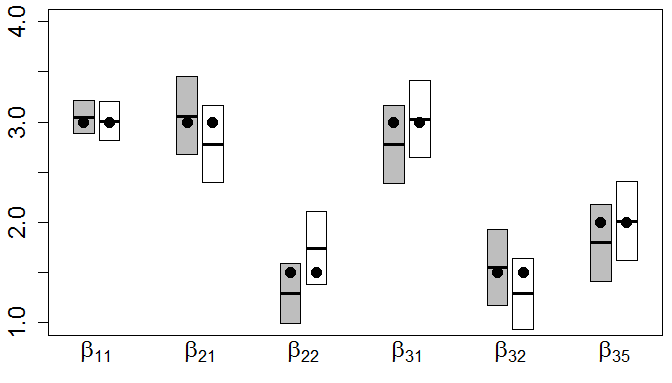}}
\subfigure[$\tau_1=100.5$, $\tau_2=150.5$, $p=500$]{\includegraphics[scale=0.3]{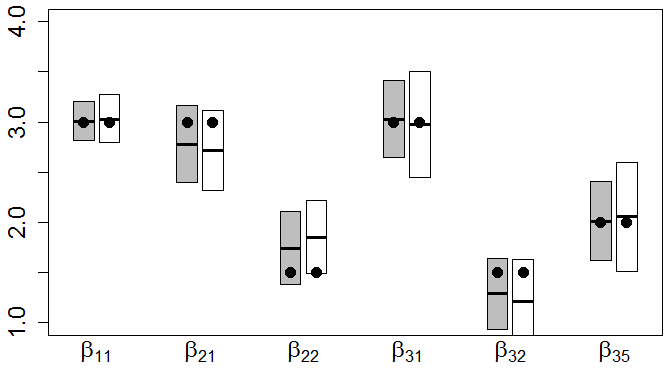}}\\
\caption{Two change point model: Posterior median estimates and $95\%$ confdence intervals of the non-zero entries of $\beta_1$, $\beta_2$ and $\beta_3$. $\beta_{kj}$ denotes the $k$ entry of $\beta_k$ for $k=1,2,3$}\label{fig:beta}
\end{figure}

Table \ref{tab:taus} presents the change-point estimates for all the scenarios. We observe that all the change points are accurately estimated.
\begin{table}[h]
\centering
\caption{Two change point model: Posterior median and $95\%$ confidence intervals of the change points.}\label{tab:taus}
\begin{tabular}{c|c|c|c|c|c} \hline
&&& (50.5,100.5) & (50.5,150.5) & (100.5,150.5) \\ \hline
\multirow{4}{*}{$p=250$} & \multirow{2}{*}{AR} & $\tau_1$ & 50.8 (49.2, 53.9) & 50.5 (48.3, 52) & 100.6 (100, 102.7)\\
 &  & $\tau_2$ & 101.5 (101, 103.3) & 148.6 (146.9, 151.7) & 150.1 (149, 152.5)\\ \cline{2-6}
 & \multirow{2}{*}{CS} & $\tau_1$ & 50.6 (45.7, 52.9) & 50.8 (49.1, 52.9) & 98.6 (92.6, 100.8)\\
 &  & $\tau_2$ & 101.3 (97.3, 104.1) & 151.3 (148.1, 152.7) & 151 (145.4, 152.8)\\ \hline
\multirow{4}{*}{$p=500$} & \multirow{2}{*}{AR} & $\tau_1$ & 50.6 (45.7, 52.9) & 50.8 (49.1, 52.9) & 98.6 (92.6, 100.8)\\
&  & $\tau_2$ & 101.3 (97.3, 104.1) & 151.3 (148.1, 152.7) & 151 (145.4, 152.8)\\ \cline{2-6}
 & \multirow{2}{*}{CS} & $\tau_1$ & 47.2 (43.3, 51.8) & 50.4 (49.1, 51.9) & 101.5 (97.6, 103.2)\\
 &  & $\tau_2$ & 102.1 (97.9, 105.8) & 147.4 (142.1, 151.8) & 146.2 (140.3, 151.3)\\ \hline
 \end{tabular}
 \end{table}
Turning to variable selection, once again, BASAD identified the exact set of predictors across all scenarios. Figure \ref{fig:beta} provides the estimates of the non-zero coefficients in the model.
We observe that with the exception of the the second entry of $\beta_2$, all the non-zero coefficients were well within the posterior confidence intervals. Overall, we observe that even for multiple change points, our methodology can accurately detect the change points, identify the correct sets of predictors, and estimate the regression coefficients.

\section{Minnesota House Price Index Data}\label{sec:minn}
In this section we apply our methodology to conduct an empirical analysis of Minnesota house price index data.

\subsection{Literature Review}\label{sec:lit}
The expansive literature on statistical analysis of house price data can be broadly classified into two major subdivisions based on their objectives. The first category of articles is aimed at forecasting individual house prices based on the constituent characteristics of the houses. Such hedonic regression models usually incorporate information regarding the structure, location, neighborhood and selling history of the house to determine its price range \citep[refer][for a comprehensive review on hedonic models]{hedonic03}. The second class of analysis focuses on 
understanding how real estate prices impact or get impacted by the economy of a country or a region. Until very recently this segment received relatively scant attention because interaction between housing market and macroeconomic variables was often deemphasized \citep{leung04}. The US sub-prime mortgage crisis between 2007 and 2009 triggered by the collapse of the housing market has recuperated interest on studying this relationship.

Several empirical analysis furnish evidence for co-movements of house price indices and other macro-economic variables like Gross Domestic Product (GDP), consumer price indices, unemployment rates, interest rates, stock price indices and so on \citep{eu03,hofmann03,tsatsaronis04,otrok05,eu13,greece15}. Multivariate regression models have been used to understand the relationship between these macro economic variables and house price index (hpi) in Ukraine \citep{kyiv}, Sweden \citep{swede} and Malaysia \citep{malay}. These analyses often assume a single underlying time-homogeneous relationship between hpi and the explanatory variables. Such an assumption may be far fetched in reality where correlation between hpi and its macroeconomic determinants may exhibit differential trends over time. For example, as noted in \cite{retro13}, the US stock market crash in the `Internet bubble burst' of 2001-2002 was not accompanied by plummeting house prices whereas in the sub-prime mortgage crisis in 2007-2009, stocks and house prices witnessed simultaneous collapse.

While the impact of macroeconomic variables on US house prices has been analysed in the literature \citep{case01,catte04} any relevant literature focusing on similar analysis at state level has eluded us.  As housing markets are local in nature \citep{local}, a state level macro-analysis may reveal trends not reflected in a similar nationwide study. The state of Minnesota is home to 18 Fortune 500 companies (\url{http://mn.gov/deed/business/locating-minnesota/companies-employers/fortune500.jsp}) and has the second highest number of Fortune 500 companies per capita. Furthermore, the Minneapolis-St. Paul metropolitan area hosts the highest number of Fortune 500 companies per capita among the 30 largest metropolitan areas in US. Hence the local industries may play a significant role in determining real estate prices in Minnesota. We use a multivariate regression model with change points to investigate the relationship between the state-level hpi of  Minnesota and both local and national macro-economic variables.

\subsection{Data and Model}\label{sec:datmod} 
\begin{figure}[h!]
\includegraphics[scale=0.28]{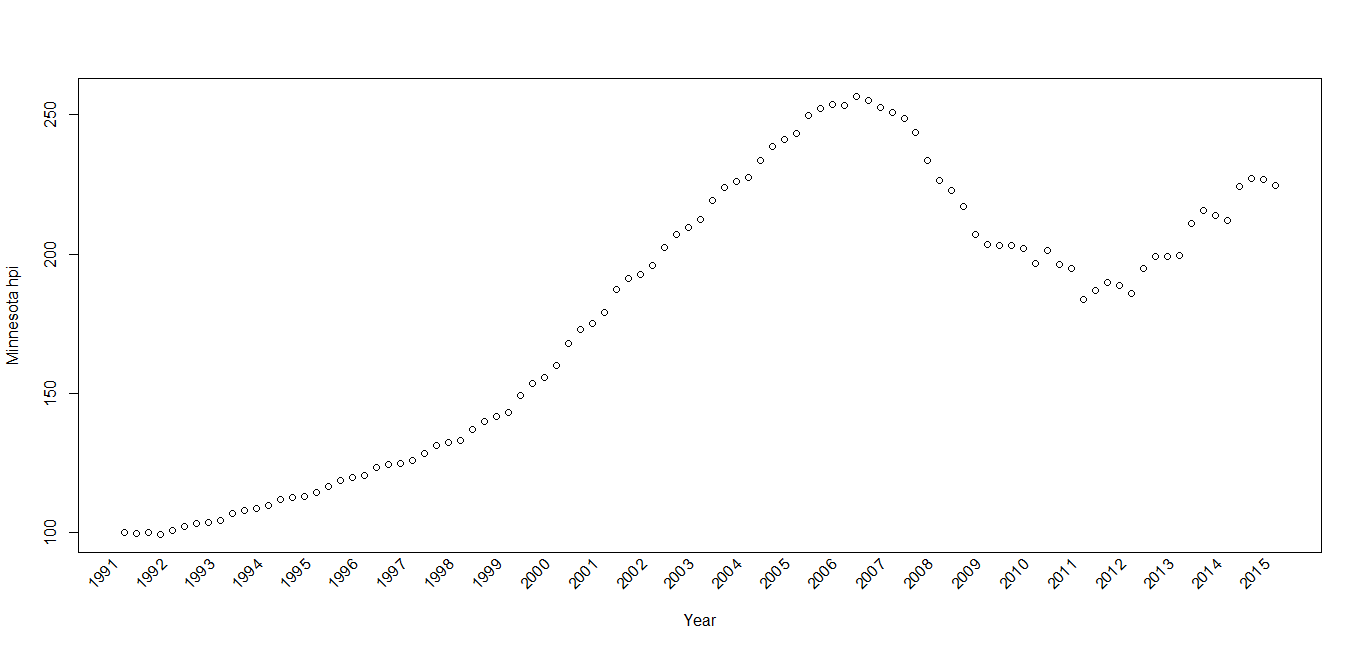}
\caption{Minnesota hpi time series}\label{fig:data}
\end{figure}

We use quarterly Minnesota hpi data published by the Federal Housing Finance Agency (FHFA) from the first quarter of 1991 to the first quarter of 2015 for our analysis. Quarterly state hpi data was obtained from \url{http://www.fhfa.gov/DataTools/Downloads/Documents/HPI/HPI_EXP_state.txt}. Figure \ref{fig:data} plots the hpi time series. We observe that there are two possible break points --- one around 2006-2008 where hpi starts to depreciate after reaching a peak and one later around 2012 where hpi starts its revival. However, this merely suggests possible change points in terms of the overall mean level for hpi. Our focus here is the relation between house price and other economic factors. We want to see if there is a change point in time such that the model before and after the change point is different.


The macro-economic indices used as explanatory variables include national unemployment rate (unemp) and national consumer price indices (cpi). The monthly cpi data was obtained from \url{http://inflationdata.com/inflation/consumer_price_index/historicalcpi.aspx?reloaded=true#Table?reloaded=true} while the unemployment data was obtained from \url{http://data.bls.gov/timeseries/LNS14000000}. All monthly indices were averaged to convert to quarterly indices.
Instead of including a national stock index  in the model like the S\&P 500 or the Dow Jones Industrial Average, we use the stock prices of Minnesota based Fortune 500 companies. 14 out of the 18 Minnesota-based Fortune 500 companies were publicly traded since before 1991 and we include their stock prices in the regression model. Additionally, the list of top 10 employers in Minnesota include Wal-Mart Stores Inc. and Wells Fargo Bank Minnesota (\url{http://inflationdata.com/inflation/consumer_price_index/historicalcpi.aspx?reloaded=true#Table}). Hence, the stock prices of these two companies are also included in the model.
The 16 stocks used in total are listed in Table \ref{tab:stocks}.

\begin{table}[h!]
\centering
\caption{List of stocks used in Minnesota hpi analysis}\label{tab:stocks}
\begin{tabular}{cc|cc} \hline
Company Name & Ticker Symbol & Company Name & Ticker Symbol \\ \hline
3M Company & MMM & St. Jude Medical, Inc. & STJ \\
Best Buy Co., Inc. & BBY & SuperValu, Inc. & SVU \\
Ecolab, Inc. & ECL & Target Corporation & TGT \\
Fastenal Co. & FAST & UnitedHealth Group Inc. & UNH \\
General Mills, Inc. & GIS & U.S. Bancorp & USB \\
Hormel Foods Corporation & HRL & Wal-Mart Stores, Inc. & WMT \\
Medtronic Plc. & MDT & Wells Fargo \& Company & WFC \\
Mosaic Company & MOS & Xcel Energy Inc. & XEL \\ \hline
\end{tabular}
\end{table}

\begin{figure}
\centering
\includegraphics[scale=0.5]{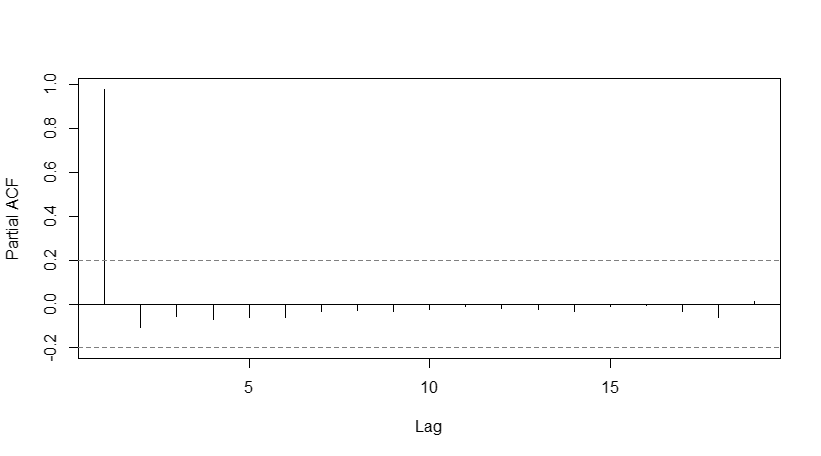}
\caption{Partial autocorrelation function for Minnesota hpi time series}\label{fig:pacf}
\end{figure}

Financial indices often exhibit strong autocorrelation and consequently autoregressive components commonly feature in house price models \citep{nagaraja2011}. Figure~\ref{fig:pacf} plots the partial auto-correlation values of the hpi time series as a function of the lag.
We observe that the index lagging one quarter behind (AR(1)) has very high correlation with the hpi time series but it quickly falls off beyond the first lag and all the subsequent lagged indices have insignificant partial correlations. Consequently, we include only the AR1 term in the regression model.

Statistical analysis involving financial time series is often preceded by customary seasonality adjustment of the indices using standard time series techniques. It is well known that house price time series reveal a predictable and repetitive pattern with systematic highs in summer and lows in winter \citep{seasons}. Consequently, publishers of  popular house price indices like the FHFA or Standard and Poor's (Case-Shiller index) produce a version of their indices discounting this effect \citep{fhfa14}. However, Minnesota is a land of extreme climates experiencing one of the widest range of temperatures in U.S. It is of interest to investigate if the impact of weather in Minnesota on its house prices extends beyond the routine pattern. Hence, we include the state level quarterly average temperatures (temp) and precipitation (precip) in the model. The data is obtained from \url{http://www.ncdc.noaa.gov/cag/time-series}.

Our model is stated as follows. Assuming $K$-change points $\tau_1 < \tau_2 < \ldots < \tau_K$ where $K$ is to de determined by the data, the regression model at the $t^{th}$ quarter, for $ \tau_{k-1} \leq t \leq \tau_k$ is given by:
\begin{align}\label{eq:model}
hpi_t&=\beta_k^{intercept}+\beta_k^{ar(1)} \mbox{hpi}_{t-1} + \beta_k^{cpi} \mbox{cpi}_t + \beta_k^{unemp} \mbox{unemp}_t \nonumber \\
& + \beta_k^{temp} \mbox{temp}_t + \beta_k^{precip} \mbox{precip}_t +  {\beta_k^{stocks}}'\mbox{stocks}_t + \eps_t
\end{align}
Here stocks$_t $ denote the $16 \times 1$ vector formed by stacking up the stock prices at time $t$ of the companies listed in Table \ref{tab:stocks} and $\beta_k^{stocks}$ is the corresponding coefficient vector.

\subsection{Results}\label{sec:results} We used the data from the second quarter of 1991 to the second quarter to 2014 for model fitting. The first quarter data of 1991 was used for the AR(1) term, whereas the data for last 2 quarters of 2014 and first quarter of 2015 were held out for out-of-sample validation. Under the assumption of $K$ change points, separate regression models are fit to each of the $K+1$ segments. Hence, although the sample size ($n =93$) was larger than the number of predictors ($p=22$) in model \ref{eq:model}, depending on the location of change-points, many segments may have less than 22 datapoints thereby necessitating high-dimensional regression methods. Higher values of $K$ ($\geq 3$) implies that average sample size for each segment ($(n/(K+1)$) becomes really small $(< 25)$ and the estimates obtained may not be reliable.  Hence, we restrict ourselves to the choices $K=0,1$ and $2$ and fit model \ref{eq:model} using the BASAD priors for the coefficient vectors in each segment. Note that for $K=0$ i.e. no change point, the model simply reduces to the traditional BASAD model. The models for different values of $K$ were assessed based on their in-sample DIC score and out-of-sample RMSPE score \citep{rmspe02}. Due to the presence of the autoregressive term, out-of-sample forecasts were obtained using one-step-ahead predictions.


Table \ref{tab:number} contains the DIC and RMSPE scores and estimated change points for the three choices of $K$. \begin{table}[h]
\centering
\caption{Minnesota hpi analysis: DIC, RMSPE scores and estimated change points }\label{tab:number}
\begin{tabular}{c|c|c|c|c} \hline
K & DIC & RMSPE & $\hat \tau_1$ & $\hat \tau_2$ \\ \hline
0 & 286.575 & 66.743 & & \\
1 & 249.587 & 14.181 & 2008Q4 (2008Q2, 2009Q1) & \\
2 & 269.054 & 14.967 & 2006Q2 (2005Q4, 2007Q3) & 2011Q1 (2010Q4, 2011Q2) \\ \hline
\end{tabular}
\end{table}
Both the DIC score and the RMSPE score for $K=0$ were significantly worse than the scores for $K=1$ and $2$ justifying the use of a change point model. The single change point model detected a change point around late 2008- early 2009 which coincides with the sub-prime mortgage crisis. The two change point model detected change points in mid 2006 and early 2011. The DIC score for the single change point model was substantially better. RMSPE scores for change point models only validate the accuracy of the models after the last change point. We observed that the RMSPE score were similar for $K=1$ and $2$ with the former turning out to be marginally better. Hence, we present the subsequent analysis only for the single change point model as both in-sample and out-of-sample validations provide strongest evidence in favor of it.

\begin{figure}[!h]
\centering
\includegraphics[scale=0.5]{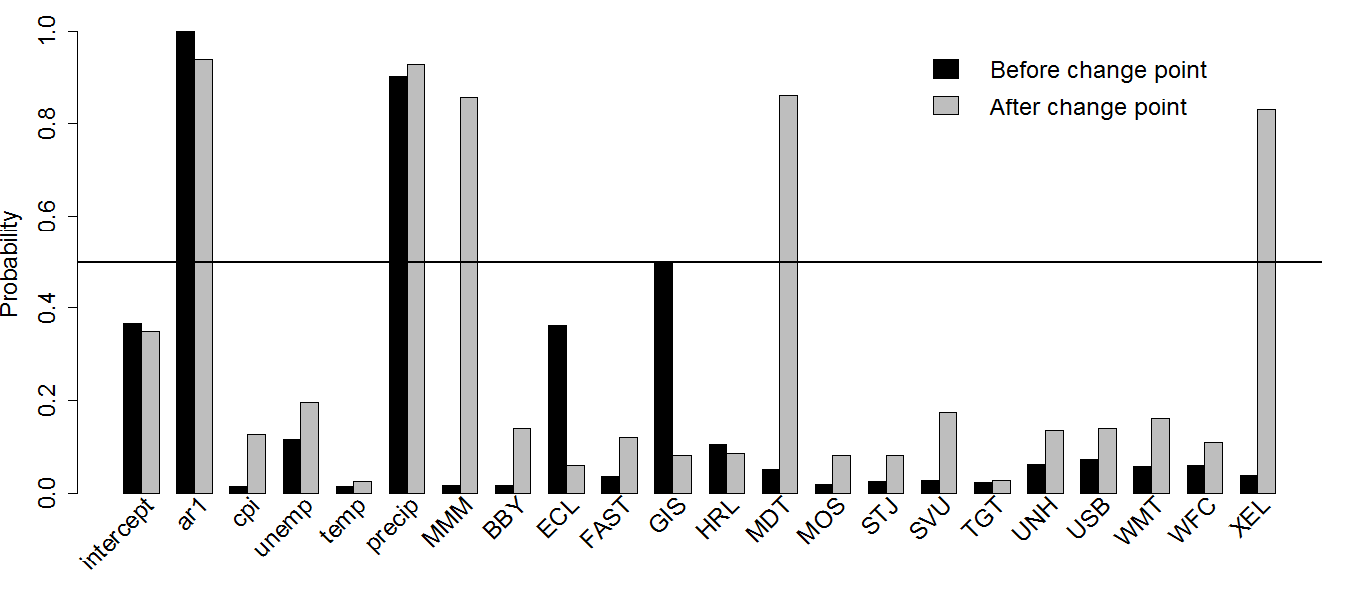}
\caption{Minnesota hpi analysis: Posterior median probabilities of variable selection using single change point model}\label{fig:prob}
\end{figure}

Figure \ref{fig:prob} plots the probability of selection each of the regressors in model (\ref{eq:model}) before and after the change point using BASAD priors.
 We observe that the set of variables selected by the median probability model differs before and after the change point. The AR(1) index and precipitation are selected with high probabilities in both segments. However, the selection of stocks differ considerably on either side of the change point. We see that prior to change point in 2008, there was little correlation between hpi and stocks with only General Mills (GIS) having a posterior median probability close to 0.5 (0.498). Perhaps this is a reflection of the fact discussed earlier that stock prices and hpi did not exhibit co-movements during the early 2000s. After the change point in 2008 the stocks of 3M (MMM), Medtronic (MDT) and Xcel Energy (XEL) are selected with high probability.

\begin{figure}[!h]
\includegraphics[width=10cm, height=4cm]{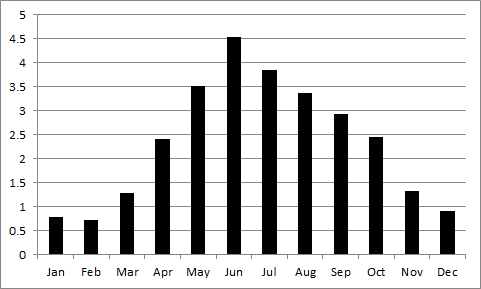}
\centering
\caption{Avg. monthly precipitation in Minnesota between 1991Q1 to 2015Q1}\label{fig:precip}
\end{figure}

Turning to the actual coefficient values presented in Table  \ref{tab:coeff} we observe that the value for coefficient for the AR(1) index drops significantly post change point indicating less autoregressive behavior after the change point.
\begin{table}[h!]
\centering
\caption{Minnesota hpi analysis: Posterior median (and 95$\%$ confidence intervals) for the coefficients. The variables which are selected in at least one segment are indicated by *.}\label{tab:coeff}
\begin{tabular}{c|c|c} \hline
& Before change point & After change point \\ \hline
intercept & -0.015 (-4.537, 2.65) & 0.004 (-3.557, 4.003) \\
AR(1)* & 1.034 (0.96, 1.106) & 0.761 (0.022, 0.873) \\
cpi & 0 (-0.055, 0.054) & 0.042 (-0.07, 0.861) \\
unemp & -0.015 (-0.697, 0.137) & 0.01 (-1.063, 1.304) \\
temp & -0.02 (-0.064, 0.044) & 0.063 (-0.021, 0.121) \\
precip* & 0.964 (-0.003, 1.63) & 1.952 (-0.012, 2.806) \\
MMM* & 0.019 (-0.065, 0.107) & -0.573 (-0.832, 0.015) \\
BBY & -0.003 (-0.091, 0.088) & 0.065 (-0.071, 0.322) \\
ECL & -0.104 (-0.541, 0.048) & 0.026 (-0.09, 0.217) \\
FAST & 0 (-0.132, 0.132) & -0.033 (-0.746, 0.093) \\
GIS* & -0.151 (-0.958, 0.068) & -0.002 (-0.206, 0.185) \\
HRL & -0.013 (-0.636, 0.131) & 0.016 (-0.12, 0.438) \\
MDT* & 0.071 (-0.022, 0.168) & 1.214 (-0.073, 1.666) \\
MOS & -0.056 (-0.107, -0.004) & -0.017 (-0.413, 0.074) \\
STJ & 0.007 (-0.097, 0.107) & 0.026 (-0.101, 0.394) \\
SVU & -0.035 (-0.144, 0.07) & 0.021 (-0.129, 1.569) \\
TGT & 0.001 (-0.108, 0.11) & 0 (-0.116, 0.117) \\
UNH & -0.061 (-0.185, 0.043) & -0.042 (-0.521, 0.083) \\
USB & 0.045 (-0.071, 0.241) & -0.017 (-2.98, 0.147) \\
WMT & 0.073 (-0.016, 0.17) & 0.014 (-0.122, 1.083) \\
WFC & -0.006 (-0.154, 0.14) & 0.001 (-0.142, 2.296) \\
XEL* & -0.005 (-0.13, 0.112) & 1.279 (-0.05, 2.436) \\ \hline
\end{tabular}
\end{table}
We also observe a positive association of hpi with precipitation. Since summer months witness significantly higher precipitation than winter (see Figure \ref{fig:precip}),
this merely corroborates the traditional 'hot season cold season' trend of house prices. What is more interesting is the fact that this effect is much more pronounced after 2008 indicating more disparity between summer and winter house prices in the post recession bear market.

\begin{figure}[!h]
\centering
\includegraphics[scale=0.5]{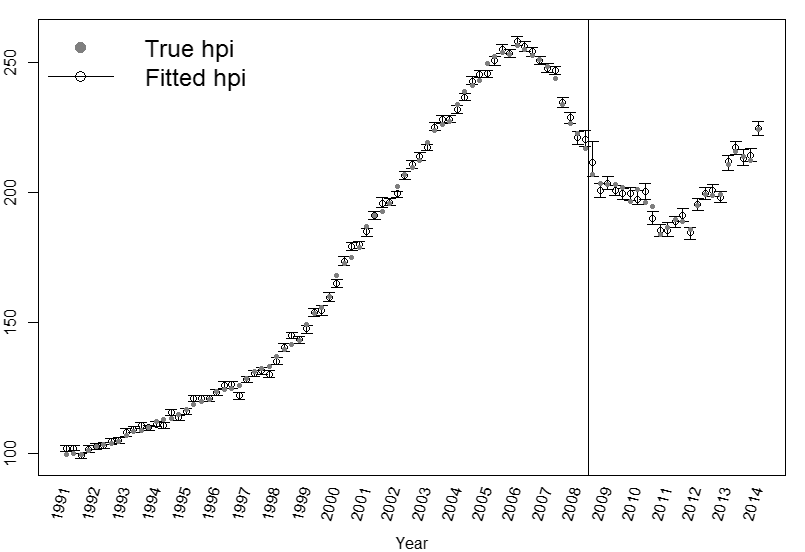}
\caption{Minnesota hpi analysis: In-sample posterior predictive medians and confidence intervals of hpi. 
 The vertical line indicates the posterior median estimate of the change point.
}\label{fig:infits}
\end{figure}

The posterior predictive distributions for the hpi at each quarter was also obtained from the MCMC sampler and Figure \ref{fig:infits} plots the true hpi versus the in-sample median fits and confidence intervals.
We see that that the posterior confidence intervals provide substantial coverage. We also observe that the confidence intervals are wider after the change point. This is expected as there are only about 23 time points compared to around 70 before the change point.
\begin{figure}[!h]
\centering
\includegraphics[scale=0.5]{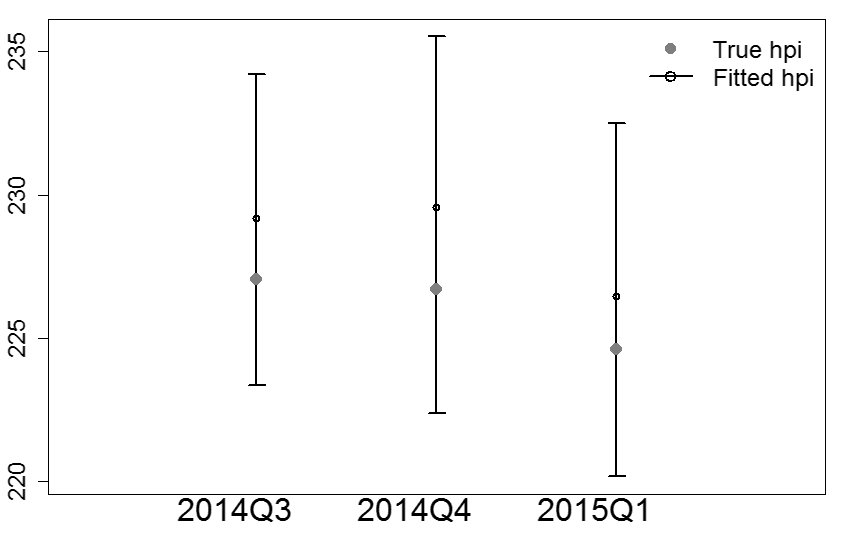}
\caption{Minnesota hpi analysis: Out-of-sample posterior predictive medians and confidence intervals of hpi. 
}\label{fig:outfits}
\end{figure}
The out-of-sample posterior fits are provided in Figure \ref{fig:outfits}. All the three out-of-sample predictions are very close to their true values. However, the out-of-sample confidence intervals are significantly wider reflecting the uncertainty associated with prediction.

Observe from Figure \ref{fig:prob} that in presence of the stock prices of Minnesota based companies, national level macro-economic indicators like the cpi or unemployment were not selected in the model. This perhaps provides evidence in support of the conjecture that hpi is strongly correlated with local macro-economics \citep{local}. However, its worthwhile to point out that multivariate regression models, although a simple and powerful tool to determine correlation, rarely implies causality. Any confirmatory assessment of the change points detected and the variables selected by our method would require further economic research. Nevertheless, the model evaluation metrics in Table \ref{tab:number} provide very strong evidence in favor of one or more change points thereby justifying the use of our methodology to analyze the data.

\section{Conclusion}\label{sec:dconc} We have presented a very general methodology for analyzing high dimensional data using changing linear regression. Our fully Bayesian approach offers several inferential advantages including quantifying uncertainty regarding the the change points as well as variable selection for each segment. Our framework is flexible to the choice of variable selection priors although the BASAD prior empirically outperformed other competing choices. A wide range of constrained variable selections like grouping or partial selection can be seamlessly accomplished in our setup. The analysis of Minnesota hpi data using our methodology revealed strong evidence for a potential change point with respect to the association with other macro-economic variables.

We have discussed several approaches for handling unknown number of change points. However, most of them comes with statutory warnings regarding computational requirements. More efficient algorithms for simultaneous detection of number of change points need to be researched.
Other potential extensions include accommodating missing data,  measurement errors or non-Gaussian responses in a high dimensional changing regression setup. Extensions to change point detection in high dimensional VAR models also need to be explored due to the extensive usage of VAR models in economics research \citep{cr00,favar}. In a time series context, our work is restricted to detecting historical change points. Detecting future change points in high dimensional time series is equally important to provide accurate predictions. We identify all these areas as directions for future research.




\bibliographystyle{asa}
\bibliography{hdcpbib}
\label{lastpage}

\end{document}